\documentclass[10pt]{article}

\usepackage{graphicx}
\usepackage[dvips]{epsfig,graphicx}
\usepackage{amsfonts,amssymb,amsmath,amsthm,multirow}%
\usepackage{pseudocode}
\usepackage{authblk}
\usepackage{algorithm}
\usepackage{algorithmic}
\usepackage{float}
\usepackage{tikz}
\usetikzlibrary{positioning,chains,fit,shapes,calc}

\rmfamily

\newtheorem{Theorem}{Theorem}[section]
\newtheorem{Definition}{Definition}[section]

\newtheorem{Proposition}[Theorem]{Proposition}
\newtheorem{Corollary}[Theorem]{Corollary}

\title {On The Relationship between (16,6,3)-Designs and (25,12) Self-Orthogonal Codes\footnote{Published in the Journal of Combinatorial Mathematics and Combinatorial Computing (JCMCC), volume 107, 137–147, 2018}}
\author[1]{Navid Nasr Esfahani}
	\author[2]{G.~H.~John van Rees}
	\affil[1]{David R.\ Cheriton School of Computer Science, University of Waterloo,
		Waterloo, Ontario N2L 3G1, Canada}
	\affil[2]{Department of Computer Science, University of Manitoba, Winnipeg, Manitoba R3T 2N2, Canada}
\date{}
\begin{document}
	\maketitle 
	\begin{abstract}
		Binary self-orthogonal codes and balanced incomplete block designs are two combinatorial configurations that have been much studied because of their wide areas of application. In this paper, we have shown the distribution of $(16,6,3)$-designs in binary (25,12) self-orthogonal codes. The paper also presents the relationships among the codes with embedded designs.
	\end{abstract}
	
	\section{Introduction}
	Balanced incomplete block designs (BIBDs) have been studied for many decades. A BIBD is a combinatorial design with $v$ elements and $b$ $k$-sets, such that each element appears in $r$ sets, and each pair of elements appear together in $\lambda$ sets. Fisher and Yates \cite{FY} listed many instances of BIBDs to be used in designing specific experiments in different fields of science. In addition to their application in statistics, BIBDs has been used in cryptography: D'Arco \textit{et al.}~\cite{DArcoAONT} used a family of BIBDs to achieve transforms close to 2-\textit{all or nothing transforms}. Formally, a BIBD is defined as the following.
	
	\begin{Definition}\label{Def_BIBD}
		A $(v,k,\lambda)$ balanced incomplete block design ($(v,k,\lambda)$-design) is a pair $(V,B)$ where $V$ is a $v$-set and $B$ is a collection of $b$ $k$-subsets, called blocks, of $V$ such that each element of $V$ is contained in exactly $r$ blocks and any 2-subset of $V$ is contained in exactly $\lambda$ blocks. 
	\end{Definition}
	
	The numbers $v$, $b$, $r$, $k$, and $\lambda$ are the parameters of the design. It can be shown that $r$ and $b$ can be calculated from the other parameters, as $rv=bk$ \cite[p.5]{S book} and $r(k-1)=\lambda(v-1)$ \cite[p.4]{S book}.
	 A $(v,k,\lambda)$-design can be represented by its \textit{incidence matrix}, i.e., a $v \times b$ matrix of binary entries. Each row of the matrix corresponds to an element of $V$, and each column represents a block. An entry in the incidence matrix is $1$ if and only if the element corresponding to its row appears in the block corresponding to its column. Therefore, there are $r$ 1's in each row and $k$ 1's in each column; also the inner product of each two rows equals $\lambda$.

	Since we will only use binary codes in this paper, we only consider binary codes and all the arithmetic on codes or designs will be \textit{modulo} $2$, unless otherwise specified. A binary code of length $n$ is a subset of $V_2(n)$, which is the set of all the binary n-tuples, called \textit{words}. The words that are in the code are called \textit{codewords}. The weight of a codeword is the number of non-zero coordinates in it. Also, the distance between two words means Hamming distance, i.e., the number of coordinates where they differ. Consequently, the distance of a code is defined to be the minimum distance between all pairs of codewords in it. A code $C$ is linear if for any two codewords $u,v\in C$, $u\oplus v\in C$, as well. Since a linear code is a subspace of a vector space, we can represent it by a set of basis that spans the subspace. A matrix is said to generate a code if its the rows of the matrix spanning the code. A \textit{generator matrix} for a code $C$ is a full rank matrix that generates $C$. An $(n,k)$ code is a linear code of length $n$ with $k$ orthogonal basis, and if we know that the code is of distance $d$, it is written as $(n,k,d)$ code.

	\begin{Definition}\label{Def_OrthComplement}
		The \textit{orthogonal complement} of a linear binary code $C$, represented by $C^{\perp}$, is the subset of $V_2(n)$ that contains all the words which are orthogonal to all codewords of $C$. 
	\end{Definition}
	
	Clearly, $C^{\perp}$ is also a binary linear code.
	
	\begin{Definition}
		A linear binary code $C$ is \textit{self-orthogonal} if $C\subseteq C^{\perp}$. The code $C$ is \textit{self-dual} if $C=C^{\perp}$.
	\end{Definition}
	
	For any code $C$ of length $n$, its \textit{weight distribution} is an array $W$ of length $n+1$, where the $i^{th}$ element of $W$ is the number of codewords of weight $i$. A code is \textit{even} if all its codewords have even weight. And finally, if all codewords in a code $C$ have weights that are multiples of $4$, then $C$ is \textit{doubly even}.

	Two binary codes are \textit{equivalent} if the coordinates of one are a permutation of the coordinates of the other one. Since reordering the coordinates of a code does not change the main characteristics of the code, we do not differentiate between equivalent codes.
	
	The study of codes and designs, and their relationship, is not a new topic. For example, Assmus and Key \cite{AK} used designs to produce block codes, that are generated by the columns of the incidence matrix. Munemasa and Tonchev \cite{MT} used the block code of a known $(56,16,6)$-design to construct a new design with the same parameters. Another method to obtain a code from a design is to use the rows of the incidence matrix to generate a point code. In particular, in this paper we will use certain designs to produce point codes that are self-orthogonal. This has been done previously in the study of the putative quasi-residual $(22,8,4)$-design, which has $\lambda=4$ in the family of $(6\lambda-2, 2\lambda,\lambda)$-designs. Bilous and van Rees \cite{BvR} proposed the idea of searching the appropriate $(33,16)$-self-orthogonal codes for the incidence matrix of a $(22,8,4)$-design.  The incidence matrix of the $(22,8,4)$-design can be used to generate a self-orthogonal doubly-even $(33,16)$ code. In the final step, Bilous \textit{et al.}~\cite{etal} used computers to search all binary self-orthogonal doubly-even $(33,16)$ codes. Since the search did not find the desired incidence matrix in any of the codes, there is no $(22,8,4)$-design.
	
	The goal of this paper is to see how $(16,6,3)$-designs are embedded in self-orthogonal doubly-even $(25,12)$ codes. To achieve that, we need the designs and codes, which are both already enumerated.
	Spence \cite{S} generously provided us with his results of generating non-isomorphic $(16,6,3)$-designs, the number of which, $18920$ non-isomorphic ones, is recorded in Handbook of Combinatorial Designs~\cite{handbook}. On the other hand, Bouyukliev kindly ran his program, Q-extension~\cite{B}, to generate self-orthogonal doubly-even $(25,12)$ codes and provide us with them. Having both of the configurations available, we classified the designs based on the rank of their incidence matrices, the weight distribution of the codes they generate, and the designs in which they are embedded. Finally, we studied the relationship of codes with embedded designs with each other.
	
	The remaining part of the paper will continue as the following. After an explanation of each step in the the procedure the results of computation will be described, and finally the distribution of designs in codes and codes with embedded designs in those of higher ranks are provided.

	\section{Designs in Codes}
	As Spence had kindly provided us with the designs, the first step was to separate the designs according to the rank of their incidence matrix, and the result of this step is given in Table~\ref{Table:design-dim}.  That the code generated by the incidence matrix of a $(16,6,3)$-design is not self-orthogonal makes it not particularly useful for our purposes. Thus, we append a column of 1's to the incidence matrix, which is then called the \textit{augmented incidence matrix}, to generate a code of length 25. The list of designs with the codes their augmented incidence matrix generate is given in \cite{BigTable}.  Since the rows of the augmented incidence matrix have even weight and their pairwise intersection is even, it is straight forward to prove that the code generated by the rows of the incidence matrix is self-orthogonal.  After generating the code this way, elementary row operations were applied to find the generator matrix of the self-orthogonal code. The rank of the generator matrix is the dimension of the code.

	\begin{table}[H]
		\begin{center}
			\caption{Distribution of Designs based on the Rank of their Incidence Matrix}
			\label{Table:design-dim}
			\begin{tabular}{|c|c|c|c|}
				\hline
				Incidence Matrix Rank & $10$ & $11$ & $12$ \\ \hline
				Number of Designs & $6$ & $245$ & $18669$ \\ \hline 
			\end{tabular}
		\end{center}
	\end{table}
	
		 By enumerating the codewords of each code from the generator matrices obtained form the augmented incidence matrices, the weight distribution of each code is computed, and the codes are classified based on their weight distribution. This primary classification of the designs is interesting by itself, but it also expedites final classification. The result of this classification is presented in Table~\ref{Table:design-weightDist}.

	However, as the following  proposition indicates, two non-equivalent codes can have the same weight distribution.
		\begin{Proposition}\label{Prop_uniqWDforEquCodes} If two codes are equivalent, they have the same weight distribution. But the same weight distribution does not imply that two codes are equivalent.
		\end{Proposition}
		\begin{proof}Changing the order of the coordinates does not change the weight of a codeword so the first part is true. In order to check the second part, one can consider the $(16,8)$ self-dual codes, generated by the generator matrices: $A_8\oplus A_8$ and $E_{16}$ using the terminology and tables from Pless \cite{P}.  Both codes share the same weight distribution array although they are not equivalent.
		\end{proof}
	
	Therefore, we needed to classify codes in each weight distribution class to classes of equivalent codes. We used \textit{nauty}, developed by McKay \cite{McKay}, which finds isomorphisms in graphs. In order to do so, each code $C$ was converted into a bipartite graph where one part of the vertices is the set of codewords and the other part is the set of coordinates positions.  Only if a codeword has a 1 in a coordinate position, the vertex representing that codeword is joined to the vertex representing that coordinate position.  In this case, the graphs will be isomorphic if and only if the corresponding codes are equivalent.

	
	\begin{table}[H]
		\begin{center}
			\caption{Distribution of Designs in Weight Distribution classes}
						\label{Table:design-weightDist}
						\setlength\tabcolsep{4pt}
						\def\arraystretch{1}
			\begin{tabular}{|c|c|c|c|c|c|c|c|c|c|c|c|c|c|c|}
			\hline
				\multicolumn{13}{|c|}{Weight Polynomial (only even weights)}& \multirow{2}{*}{$\#$} & \multirow{2}{*}{k}  \\ \cline{1-13}
				\textit{No.}&$0$&$ 2$&$ 4$&$ 6$&$ 8$&$ 10$&$ 12$&$ 14$&$ 16$&$ 18$&$ 20$&$ 22$& & \\ \hline
				$1$&$1$&$ 0$&$ 1$&$ 4$&$ 82$&$ 164$&$ 346$&$ 300$&$ 77$&$ 44$&$ 5$&$ 0$&  $1$ & $10$\\ \hline
				$2$&$1$&$ 0$&$ 3$&$ 3$&$ 78$&$ 166$&$ 346$&$ 300$&$ 81$&$ 42$&$ 3$&$ 1$&  $2$ & $10$\\ \hline
				$3$&$1$&$ 0$&$ 0$&$ 3$&$ 90$&$ 166$&$ 328$&$ 300$&$ 93$&$ 42$&$ 0$&$ 1$&  $3$ & $10$\\ \hline
				$4$&$1$&$ 0$&$ 1$&$ 9$&$ 182$&$ 324$&$ 654$&$ 606$&$ 185$&$ 84$&$ 1$&$ 1$& $4$ & $11$\\ \hline
				$5$&$1$&$ 0$&$ 1$&$ 6$&$ 180$&$ 334$&$ 660$&$ 594$&$ 179$&$ 90$&$ 3$&$ 0 $ &$8$ & $11$\\ \hline
				$6$&$1$&$ 0$&$ 2$&$ 9$&$ 178$&$ 324$&$ 660$&$ 606$&$ 181$&$ 84$&$ 2$&$ 1$& $8$ & $11$\\ \hline
				$7$&$1$&$ 0$&$ 1$&$ 11$&$ 180$&$ 318$&$ 660$&$ 612$&$ 179$&$ 82$&$ 3$&$ 1$& $8$ & $11$\\ \hline
				$8$&$1$&$ 0$&$ 3$&$ 9$&$ 174$&$ 324$&$ 666$&$ 606$&$ 177$&$ 84$&$ 3$&$ 1$& $11$ & $11$\\ \hline
				$9$&$1$&$ 0$&$ 2$&$ 6$&$ 176$&$ 334$&$ 666$&$ 594$&$ 175$&$ 90$&$ 4$&$ 0$& $12$ & $11$\\ \hline
				$10$&$1$&$ 0$&$ 0$&$ 8$&$ 182$&$ 328$&$ 660$&$ 600$&$ 177$&$ 88$&$ 4$&$ 0$& $12$ & $11$\\ \hline
				$11$&$1$&$ 0$&$ 2$&$ 11$&$ 176$&$ 318$&$ 666$&$ 612$&$ 175$&$ 82$&$ 4$&$ 1$& $12$ & $11$\\ \hline
				$12$&$1$&$ 0$&$ 7$&$ 8$&$ 154$&$ 328$&$ 702$&$ 600$&$ 149$&$ 88$&$ 11$&$ 0$& $12$ & $11$\\ \hline
				$13$&$1$&$ 0$&$ 2$&$ 8$&$ 174$&$ 328$&$ 672$&$ 600$&$ 169$&$ 88$&$ 6$&$ 0$& $16$ & $11$\\ \hline
				$14$&$1$&$ 0$&$ 4$&$ 6$&$ 168$&$ 334$&$ 678$&$ 594$&$ 167$&$ 90$&$ 6$&$ 0$& $16$ & $11$\\ \hline
				$15$&$1$&$ 0$&$ 4$&$ 11$&$ 168$&$ 318$&$ 678$&$ 612$&$ 167$&$ 82$&$ 6$&$ 1$& $16$ & $11$\\ \hline
				$16$&$1$&$ 0$&$ 5$&$ 9$&$ 166$&$ 324$&$ 678$&$ 606$&$ 169$&$ 84$&$ 5$&$ 1$& $20$ & $11$\\ \hline
				$17$&$1$&$ 0$&$ 1$&$ 8$&$ 178$&$ 328$&$ 666$&$ 600$&$ 173$&$ 88$&$ 5$&$ 0$& $26$ & $11$\\ \hline
				$18$&$1$&$ 0$&$ 3$&$ 8$&$ 170$&$ 328$&$ 678$&$ 600$&$ 165$&$ 88$&$ 7$&$ 0$& $64$ & $11$\\ \hline
				$19$&$1$&$ 0$&$ 39$&$ 36$&$ 234$&$ 592$&$ 1486$&$ 1272$&$ 261$&$ 144$&$ 27$&$ 4$& $4$ & $12$\\ \hline
				$20$&$1$&$ 0$&$ 21$&$ 36$&$ 306$&$ 592$&$ 1378$&$ 1272$&$ 333$&$ 144$&$ 9$&$ 4$& $8$ & $12$\\ \hline
				$21$&$1$&$ 0$&$ 15$&$ 36$&$ 330$&$ 592$&$ 1342$&$ 1272$&$ 357$&$ 144$&$ 3$&$ 4$& $16$ & $12$ \\ \hline
				$22$&$1$&$ 0$&$ 12$&$ 36$&$ 342$&$ 592$&$ 1324$&$ 1272$&$ 369$&$ 144$&$ 0$&$ 4$& $26$ & $12$\\ \hline
				$23$&$1$&$ 0$&$ 0$&$ 21$&$ 378$&$ 640$&$ 1288$&$ 1218$&$ 381$&$ 168$&$ 0$&$ 1$& $30$& $12$\\ \hline
				$24$&$1$&$ 0$&$ 17$&$ 26$&$ 314$&$ 624$&$ 1378$&$ 1236$&$ 325$&$ 160$&$ 13$&$ 2$& $115$ & $12$\\ \hline
				$25$&$1$&$ 0$&$ 14$&$ 26$&$ 326$&$ 624$&$ 1360$&$ 1236$&$ 337$&$ 160$&$ 10$&$ 2$& $120$ & $12$\\ \hline
				$26$&$1$&$ 0$&$ 19$&$ 16$&$ 298$&$ 656$&$ 1414$&$ 1200$&$ 293$&$ 176$&$ 23$&$ 0$& $172$ & $12$\\ \hline
				$27$&$1$&$ 0$&$ 5$&$ 26$&$ 362$&$ 624$&$ 1306$&$ 1236$&$ 373$&$ 160$&$ 1$&$ 2$& $216$ & $12$\\ \hline
				$28$&$1$&$ 0$&$ 13$&$ 16$&$ 322$&$ 656$&$ 1378$&$ 1200$&$ 317$&$ 176$&$ 17$&$ 0$& $292$ & $12$\\ \hline
				$29$&$1$&$ 0$&$ 11$&$ 26$&$ 338$&$ 624$&$ 1342$&$ 1236$&$ 349$&$ 160$&$ 7$&$ 2$& $612$ & $12$\\ \hline
				$30$&$1$&$ 0$&$ 12$&$ 21$&$ 330$&$ 640$&$ 1360$&$ 1218$&$ 333$&$ 168$&$ 12$&$ 1$& $624$ & $12$\\ \hline
				$31$&$1$&$ 0$&$ 1$&$ 16$&$ 370$&$ 656$&$ 1306$&$ 1200$&$ 365$&$ 176$&$ 5$&$ 0$& $836$ & $12$\\ \hline
				$32$&$1$&$ 0$&$ 9$&$ 21$&$ 342$&$ 640$&$ 1342$&$ 1218$&$ 345$&$ 168$&$ 9$&$ 1$& $1032$ & $12$\\ \hline
				$33$&$1$&$ 0$&$ 8$&$ 26$&$ 350$&$ 624$&$ 1324$&$ 1236$&$ 361$&$ 160$&$ 4$&$ 2$& $1086$ & $12$\\ \hline
				$34$&$1$&$ 0$&$ 10$&$ 16$&$ 334$&$ 656$&$ 1360$&$ 1200$&$ 329$&$ 176$&$ 14$&$ 0$& $1408$ & $12$\\ \hline
				$35$&$1$&$ 0$&$ 3$&$ 21$&$ 366$&$ 640$&$ 1306$&$ 1218$&$ 369$&$ 168$&$ 3$&$ 1$& $1818$ & $12$\\ \hline
				$36$&$1$&$ 0$&$ 7$&$ 16$&$ 346$&$ 656$&$ 1342$&$ 1200$&$ 341$&$ 176$&$ 11$&$ 0$& $2106$ & $12$\\ \hline
				$37$&$1$&$ 0$&$ 6$&$ 21$&$ 354$&$ 640$&$ 1324$&$ 1218$&$ 357$&$ 168$&$ 6$&$ 1$& $3814$ & $12$\\ \hline
				$38$&$1$&$ 0$&$ 4$&$ 16$&$ 358$&$ 656$&$ 1324$&$ 1200$&$ 353$&$ 176$&$ 8$&$ 0$& $4334$ & $12$\\ \hline
			\end{tabular}
		\end{center}
	\end{table}

	After converting the codes into graphs, the corresponding graphs were compared and classified using {\tt shortg} procedure in \textit{nauty}. The designs generate 44 non-equivalent codes. The generator matrix, number designs generating the code, the weight polynomial number, and some other pieces of information regarding the codes are provided in \cite{NNEMT}.
	In the next step, the codes of dimensions $10$ or $11$ are embedded in codes of dimension $12$ so that the $(16,9,3)$-designs can all be found in the dimension 12 codes. However, it is unnecessary to consider all binary self-orthogonal $(25,12)$ codes, as some of them cannot contain these designs as the following theorem by Mann~\cite{Mann} shows.
	
	\begin{Theorem}
		Let $D$ be a balanced incomplete block design with parameters $v$, $b$, $r$, $k$, and $\lambda$. If $s$ blocks of $D$ are identical and
		$r > \lambda$, then $\frac{r}{k} = \frac{b}{v} \ge s$.
	\end{Theorem}
	
	Hence we immediately get the following useful corollary,
	\begin{Corollary}\label{lemma}
		There are no identical blocks in a $(6\lambda-2, 9\lambda-3, 3\lambda, 2\lambda,\lambda)$-design.
	\end{Corollary}

	\begin{Theorem}
		If a binary self-orthogonal $(25,12)$ code contains $16$ codewords that form an augmented incidence matrix of a $(16,6,3)$-design, then the code has distance greater than or equal to $4$ and has no all-0-coordinate.
	\end{Theorem}
	
	\begin{proof}
		The incidence matrix has 1's in it and so does the column that was attached so the code cannot have an all-0-coordinate. If the distance of the code is less than 4 then it must be 2 as the codewords are orthogonal.  Then the coordinate positions containing the 1's of a weight 2 codeword in the incidence matrix must always be the same.  This implies that the BIBD has 2 identical blocks.  This contradicts Corollary~\ref{lemma}.
	\end{proof}
	
	These self-orthogonal $(25,12)$ codes with distance greater than 2 with no all-0-coordinate codes are relatively easy to produce from a list of all self-orthogonal $(25,12)$ codes or more efficiently, can be produced directly.  We can now state the final result of our programs for this section.

		After embedding the codes of dimension less than $12$ with designs in codes of dimension 12, there are 36 non-equivalent binary self-orthogonal $(25,12,\geq 4)$ codes with no all-0-coordinate containing 16 rows that form the augmented incidence matrix of a $(16,6,3)$-design. The 36 codes are made up of the original 20 codes of dimension 12 and 16 new codes which contain the designs that generated the dimension 10 and 11 codes. So each code of dimension 12 with embedded designs contains on average $19348/36$ or approximately 537 non-isomorphic designs, from a minimum of 3 to a maximum of 4470.  Note that the number is 19348 as the designs whose incidence matrices are rank 10 or rank 11 can be in several codes. See Table \ref{TAB_UpdateDesNum} for the details of codes and Table \ref{TAB_CodeDist} for information about how the codes of lower dimensions with embedded designs sit in codes of dimension 12.
	
	\begin{table}[H]
		\begin{center}
			\caption{Number of Non-Isomorphic $(16,6,3)$-Designs in $(25,12)$ Self-Orthogonal Codes}
			\label{TAB_UpdateDesNum}
			\begin{tabular}{|c|c|c||c|c|c||c|c|c|}
				\hline
				Code & $k$  & Des. & Code & $k$ & Des. & Code & $k$  & Des.\\ \hline
				$C_{1}$ 	&     10 	& 1  & $C_{25}$ 	&     12 	& 4     & $C_{49}$ 	&     11 	& 2      \\ \hline
				$C_{2}$ 	&     10 	& 2  & $C_{26}$ 	&     12 	& 8     & $C_{50}$ 	&     11 	& 2  \\ \hline
				$C_{3}$ 	&     10 	& 2  & $C_{27}$ 	&     12 	& 16    & $C_{51}$ 	&     11 	& 2  \\ \hline
				$C_{4}$ 	&     10 	& 1  & $C_{28}$ 	&     12 	& 26    & $C_{52}$ 	&     11 	& 1  \\ \hline
				$C_{5}$ 	&     11 	& 6  & $C_{29}$ 	&     12 	& 31    & $C_{53}$ 	&     11 	& 1  \\ \hline
				$C_{6}$ 	&     11 	& 8  & $C_{30}$ 	&     12 	& 115   & $C_{54}$ 	&     11 	& 1  \\ \hline
				$C_{7}$ 	&     11 	& 8  & $C_{31}$ 	&     12 	& 120   & $C_{55}$ 	&     12 	& 9  \\ \hline
				$C_{8}$ 	&     11 	& 8  & $C_{32}$ 	&     12 	& 184   & $C_{56}$ 	&     12 	& 18  \\ \hline
				$C_{9}$ 	&     11 	& 8  & $C_{33}$ 	&     12 	& 216   & $C_{57}$ 	&     12 	& 10  \\ \hline
				$C_{10}$ 	&     11 	& 4  & $C_{34}$ 	&     12 	& 292   & $C_{58}$ 	&     12 	& 24  \\ \hline
				$C_{11}$ 	&     11 	& 12 & $C_{35}$ 	&     12 	& 612   & $C_{59}$ 	&     12 	& 23  \\ \hline
				$C_{12}$ 	&     11 	& 12 & $C_{36}$ 	&     12 	& 652   & $C_{60}$ 	&     12 	& 3  \\ \hline
				$C_{13}$ 	&     11 	& 12 & $C_{37}$ 	&     12 	& 862   & $C_{61}$ 	&     12 	& 24  \\ \hline
				$C_{14}$ 	&     11 	& 12 & $C_{38}$ 	&     12 	& 1050  & $C_{62}$ 	&     12 	& 12  \\ \hline
				$C_{15}$ 	&     11 	& 16 & $C_{39}$ 	&     12 	& 1086  & $C_{63}$ 	&     12 	& 4  \\ \hline
				$C_{16}$ 	&     11 	& 16 & $C_{40}$ 	&     12 	& 1476  & $C_{64}$ 	&     12 	& 8  \\ \hline
				$C_{17}$ 	&     11 	& 16 & $C_{41}$ 	&     12 	& 1848  & $C_{65}$ 	&     12 	& 52  \\ \hline
				$C_{18}$ 	&     11 	& 12 & $C_{42}$ 	&     12 	& 2153  & $C_{66}$ 	&     12 	& 16  \\ \hline
				$C_{19}$ 	&     11 	& 10 & $C_{43}$ 	&     12 	& 3874  & $C_{67}$ 	&     12 	& 20  \\ \hline
				$C_{20}$ 	&     11 	& 20 & $C_{44}$ 	&     12 	& 4470  & $C_{68}$ 	&     12 	& 6  \\ \hline
				$C_{21}$ 	&     11 	& 6  & $C_{45}$ 	&     11 	& 1     & $C_{69}$ 	&     12 	& 12  \\ \hline
				$C_{22}$ 	&     11 	& 3	 & $C_{46}$ 	&     11 	& 1     & $C_{70}$ 	&     12 	& 12  \\ \hline
				$C_{23}$ 	&     11 	& 52 & $C_{47}$ 	&     11 	& 1     & & & \\ \hline
				$C_{24}$ 	&     11 	& 11 & $C_{48}$ 	&     11 	& 2     & & & \\ \hline

			\end{tabular}
		\end{center}
	\end{table}

Finally we need to take into account the codes that contain no designs. To find the number of such codes, we need to find the number of binary $(25,12)$ self-orthogonal codes or better yet the number of binary $(25,12,\geq 4)$ self-orthogonal codes with no all-0-coordinate. Bouyukliev generated those codes using his program, Q-extension~\cite{B}, to find this information. We also computed the number of the desirable codes by using the tables of the self-dual codes by Harada and Munemasa~\cite{HM}. We cross-sectioned each $(26,13)$ self-dual code at each coordinate position to get the $(25,12)$ self-orthogonal codes and from there the non-isomorphic binary $(25,12,\geq 4)$ self-orthogonal codes with no all-0-coordinate. To cross-section a linear code at a coordinate $j$, we delete all codewords with a 1 in that coordinate and then delete the coordinate. The procedure can be reversed by introducing a new coordinate position and then adding the all one's codeword to every old codeword to get the new codewords. The new and old codewords comprise the new code.	The results of the two approaches agree and is shown in Table~\ref{Table:finaltable}.
\begin{table}[H]
	\caption{Distribution of Codes of lower dimensions with embedded designs, in codes of dimension 12}
	\label{TAB_CodeDist}
	\begin{center}
		\begin{tabular}{|c|c|c||c|c|c|}\hline
			Codes & Sub-code & $k$ & Codes & Sub-code & $k$ \\ \hline
			$C_{22},C_{24},C_{45},C_{46},C_{47}$ & $C_{1}$ & $10$ & $C_{62},C_{36},C_{43}$& $C_{18}$ & $ 11$ \\ \hline
			$C_{19},C_{48},C_{49}$ & $C_{2}$ & $10$  & $C_{57}$& $C_{19}$ & $ 11$ \\ \hline
			$C_{5},C_{50},C_{51}$ & $C_{3}$ & $10$ & $C_{67},C_{44}$& $C_{20}$ & $ 11$ \\ \hline
			$C_{9},C_{52},C_{53},C_{54}$ & $C_{4}$ & $10$ & $C_{68},C_{42},C_{37}$& $C_{21}$ & $ 11$ \\ \hline
			$C_{56},C_{41}$& $C_{5} $ & $ 11$ & $C_{60},C_{42}$& $C_{22}$ & $ 11$ \\ \hline
			$C_{61},C_{44},C_{37}$& $C_{6} $ & $ 11$ & $C_{65},C_{40},C_{44}$& $C_{23}$ & $ 11$ \\ \hline
			$C_{64},C_{43}$& $C_{7} $ & $ 11$ & $C_{59},C_{42}$& $C_{24}$ & $ 11$ \\ \hline
			$C_{55},C_{43},C_{41}$& $C_{8} $ & $ 11$ & $C_{59}$& $C_{45}$ & $ 11$ \\ \hline
			$C_{58}$& $C_{9} $ & $ 11$ & $C_{59},C_{60}$& $C_{46}$ & $ 11$ \\ \hline
			$C_{63},C_{43}$& $C_{10}$ & $ 11$ & $C_{60}$& $C_{47}$ & $ 11$ \\ \hline
			$C_{59},C_{44}$& $C_{11}$ & $ 11$ & $C_{57}$& $C_{48}$ & $ 11$ \\ \hline
			$C_{70},C_{44},C_{37}$& $C_{12}$ & $ 11$ & $C_{57}$& $C_{49}$ & $ 11$ \\ \hline
			$C_{56},C_{43}$& $C_{13}$ & $ 11$ & $C_{56}$& $C_{50}$ & $ 11$ \\ \hline
			$C_{69},C_{42},C_{32}$& $C_{14}$ & $ 11$ & $C_{56}$& $C_{51}$ & $ 11$ \\ \hline
			$C_{61},C_{42},C_{44}$& $C_{15}$ & $ 11$ & $C_{55},C_{58}$& $C_{52}$ & $ 11$ \\ \hline
			$C_{66},C_{40},C_{44}$& $C_{16}$ & $ 11$ & $C_{55},C_{29}$& $C_{53}$ & $ 11$ \\ \hline
			$C_{58},C_{36},C_{43}$& $C_{17}$ & $ 11$ & $C_{58}$& $C_{54}$ & $ 11$ \\ \hline			
		\end{tabular}
	\end{center}
\end{table}

	\begin{table}[H]
		\begin{center}
			\caption{(25,12) Self-orthogonal Codes}
			\label{Table:finaltable}
			\begin{tabular}{|c|c|c|c|c|c|}
				\hline
				Distance & 2 & 4 & 6 & 8 & total \\ \hline
				With 0 coordinate & 25 & 28 & 1 & 1 & 55 \\ \hline
				without 0 coordinate & 105 & 168 & 3 & 0 & 276 \\ \hline
				Total & 130 & 196 & 4 & 1 & 331 \\ \hline
			\end{tabular}
		\end{center}
	\end{table}

	Therefore, we can say that there are 331 inequivalent binary $(25,12)$ self-orthogonal codes, 171 of which are $(25,12,\geq 4)$ codes, and 36 of the 171 codes of the latter group contain at least one $(16,6,3)$-design.

	\section{Analysis and Conclusion}
	In this article, we have found which self-orthogonal codes contain which $(16,6,3)$-designs.  Before  we continue by discussing the larger designs in the $(6\lambda-2, 2\lambda, \lambda)$ family, let us examine this family of designs for smaller $\lambda$, i.e., $\lambda=1$ and 2.
	
	There is only one non-isomorphic $(4,2,1)$-design which is a residual design.
	Since $\lambda=1$ is odd, a column of 1's should be added to the incidence matrix whose rows will then generate a $(7,3)$ self-orthogonal code. Pless \cite{P} showed that there are two possible such codes.  Since the augmented incidence matrix has 4 rows of weight 4, the containing code must be the code with weight distribution W[0]=1 and W[4]=7.
	
	For the case $\lambda=2$, there are three non-isomorphic $(10,4,2)$-designs (each one a residual design), whose incidence matrices generate inequivalent $(15,7)$ self-orthogonal codes. These three codes have dimensions 5, 6 and 7 and the codes of dimension 5 and 6 are embedded in the code of dimension 7.  However, there are 10 inequivalent $(15,7)$ self-orthogonal codes as shown by Pless \cite{P}.  But only 4 of these self-orthogonal codes have distance greater than 2 with no all-0-coordinate. Two of these can be eliminated because they have no codewords of weight 6 where as the design has an incidence matrix whose rows all have weight 6.  Hence there are only 2 $(15,7,\geq 4)$ self-orthogonal codes with no all-0-coordinate that may contain $(10,4,2)$-designs and hence it contains all 3 non-isomorphic $(10,4,2)$-designs. If $\lambda=1$ or 2, then one self-orthogonal codes with minimum distance 4 and with no all-0-coordinate contains all the non-isomorphic designs.
	
	For the case of $\lambda=3$, each $(25,12)$ self-orthogonal code which contains a $(16,6,3)$-design has, on average, 537 non-isomorphic designs in it, with a minimum of 3 designs and a maximum of 4470 designs.  This may be an indication of a method of obtaining new designs from old designs; i.e., simply by generating the point codes, selecting the codewords of the right weight and searching those codewords for more embedded designs.
	
	For example, if the first $(16,6,3)$-design found, which was found by Bhattacharya~\cite{Bh}, was the only $(16,6,3)$-design known, then computing the augmented point code would give $C_{44}$ as listed in \cite{NNEMT}.  But that code has 4469 more non-isomorphic designs in it. But $C_6$, $C_{11}$, $C_{12}$, $C_{15}$, $C_{16}$, $C_{20}$ and $C_{23}$ are embedded in $C_{44}$. But the codes of dimension 11 are embedded in other codes such as $C_{37}$, $C_{40}$ and $C_{42}$ with another 842, 1408, and 2137 new designs respectively. Repeating this process on these codes will lead to another 172 designs in $C_{32}$.  So from 1 design, 9,029 non-isomorphic designs are obtained, just under half of all possible.
	
	Based on the distribution of designs of the $(6\lambda-2,2\lambda,\lambda)$ family in codes for the cases where the distribution of designs in codes is known, $\lambda=1,2,3,4$, the next case to be studied is $(28,10,5)$-designs. The existence of these designs is known as van Lint \textit{et al.}~\cite{vanLint} constructed 3 non-isomorphic instances of such designs. One possible extension of this work will be to generated the point codes of those designs and search them for more designs.
	
	\bibliographystyle{apsrev}

\end{document}